\documentclass[a4paper,12pt]{article}
\usepackage{amsmath,amsthm,amsfonts,amssymb,color}
\usepackage[matrix,arrow]{xy}
\usepackage{graphicx}
\usepackage[UKenglish]{babel}
\usepackage{multicol} 
\usepackage{appendix}
\usepackage{epstopdf}

\usepackage{caption}
\usepackage{subcaption}

\usepackage{multirow} 
\usepackage{float}
\usepackage{graphicx}
\usepackage{longtable}  

\makeatletter
\renewcommand\@biblabel[1]{#1.}
\makeatother

\usepackage{amsmath,amssymb,amsthm}
\usepackage{enumerate}
\usepackage{amsfonts}
\usepackage{graphics,graphicx}
\usepackage{verbatim}
\usepackage{latexsym,makeidx}
\usepackage{amsmath,amscd}
\usepackage{appendix}
\usepackage[latin1]{inputenc}

\newenvironment{keywords}{
  \vspace{2mm}
  \noindent
  \keywordsname: 
  \itshape\small
}

\def\keywordsname{\textbf{Keywords}}
  \def\mathsubclassname{\textbf{2010 AMS Subject Classification}}

\newenvironment{Figure}
  {\smallskip\noindent\minipage{\linewidth}}
  {\endminipage\smallskip}

\newtheorem{teo}{Theorem}[section]

\newcommand\be{\begin{equation}}
\newcommand\ee{\end{equation}}

\begin{document}

\title{One-phase Stefan problem with a latent heat depending on the position of the free boundary and its rate of change}
\author{
Julieta Bollati$^{1}$,  Domingo A. Tarzia $^{1}$\\ \\
\small {{$^1$} Depto. Matem\'atica - CONICET, FCE, Univ. Austral, Paraguay 1950} \\  
\small {S2000FZF Rosario, Argentina.}\\
\small{Email: JBollati@austral.edu.ar; DTarzia@austral.edu.ar.} 
}
\date{}

\maketitle
\abstract{
From the one-dimensional consolidation of fine-grained soils with threshold gradient, it can be derived a special type of Stefan problems where the seepage front, due to the presence of  this threshold gradient, exhibits the features of a moving boundary. In this kind of problems, in contrast with the classical Stefan problem, the latent heat is considered to depend inversely with the rate of change of the seepage front. In this paper a one-phase Stefan problem with a latent heat that not only depends on the rate of change of the free boundary but also on its position is studied. The aim of this analysis is to extend prior results, finding an analytical solution that recovers, by specifying some parameters, the solutions that have already been examined in the literature regarding Stefan problems with variable latent heat. Computational examples will be presented in order to examine the effect of this parameters on the free boundary. 
}


\begin{scriptsize}

\begin{longtable}{|p{8mm} p{0.5mm} p{142mm}| } 
\hline 
& &\\
{\large \textbf{Nomenclature}}&   &\\
& &\\
$a^2$  & & Diffusivity coefficient ($a^2=k/\rho c$), $[m^2/s]$.\\
$c$			\qquad	& &  Specific heat capacity, $[m^2/^{\circ}C s^2]$\\
$k$            & & Thermal conductivity, $[W/(m ^{\circ}C)]$. \\
$s$ &  & Position of the free front, $[m]$.\\
$t$& & Time, $[s]$.\\
$T$ & & Temperature, $[^{\circ}C]$.\\
$T_0$ & & Coefficient that characterizes the temperature at the fixed face, $[^{\circ}C/s^{\alpha/2}]$.\\
$x$ & & Spatial coordinate, $[m]$.\\
& &\\
$\text{Greek symbols} $ &   &\\
& &\\
$\alpha$ & & Power of the time that characterizes the temperature at the fixed boundary, dimensionless.\\
$\beta$ & & Power of the position that characterizes the latent heat per unit volume, dimensionless.\\
$\delta$ & & Power of  the velocity that characterizes the latent heat per unit volume, dimensionless.\\
$\gamma$ & & Coefficient that characterizes the latent heat per unit volume, $[s^{\delta-2} kg /(m^{\beta+\delta+1})]$.\\
$\rho$ & & Density, $[kg/ m^3]$. \\
$\xi$ & & Coefficient that characterizes the free interface, dimensionless.\\
$\eta$ & & Similarity variable in expression (\ref{Transform-1}), dimensionless.\\
& &\\
\hline
\end{longtable}

\end{scriptsize}

\newpage

\section{Introduction.}
$ $

This paper follows on from the work done by Zhou et al. in \cite{ZBL}, where a  one-dimensional consolidation process with a threshold gradient is studied. This problem is essentially a moving boundary problem where the seepage front, which moves downward gradually, plays the role of the free boundary due to the presence of this threshold gradient. This kind of problems are known in the literature as Stefan problems. They have been widely studied in the last century due to the fact that they arise in many significant areas of engineering, geoscience and industry (\cite{AlSo}-\cite{Gu},\cite{Lu},\cite{Ru},\cite{Ta4}). A review of the literature on this topic  was presentend in \cite{Ta2}.

The classical Stefan problem intends to describe the process of a material undergoing a phase change like, for example, the melting process on an ice bar.  Finding a solution to this problem consists in solving the heat-conduction equation in an unknown region  which has also to be determined, imposing an initial condition, boundary conditions and the Stefan condition at the interface.

In this paper, it will be considered the one-phase Stefan problem in a semi-infinite material with variable latent heat. The reduction to a one-phase problem is referred to the case in which it is assumed that one of the phases is at the phase-change temperature. The new mathematical feature of the problem to be solved is concerned with  the fact that the latent heat is assumed to depend on the position of the free boundary as well as on its rate of change which is our novelty, i.e $L=L(s(t),\dot{s}(t))$ (where $s(t)$ is the free front).
It is known that in the classical formulation of the Stefan problem the latent heat is constant. The idea of a variable latent heat is motivated by the following previous works:

\begin{itemize}

\item  In  \cite{Pr} it was considered a Stefan problem with a latent heat given as a function of the position of the interface $L=\varphi(s(t))$. Such assumption corresponds to the practical case when the influence of phenomena such as surface tension, pressure gradients and nonhomogenity of materials are taken into account. Sufficient conditions that ensure the existence and uniqueness of solution were studied in this paper.

\item  In \cite{VSP} it was provided an analytical solution to the one-phase Stefan problem with a latent heat defined as a linear function of the position, i.e $L=\gamma s(t)$ (with $\gamma$  a given constant). This hypothesis  makes  physical sense in the study of shoreline movement in a sedimentary basin. 
The extension to the two-phase problem was done in \cite{SaTa} .

\item  In \cite{ZWB} it was considered a one-phase Stefan problem where the latent heat is not constant but, rather a power function of the position, i.e. $L=\gamma s^{n}(t)$ (with $\gamma$ a given constant and $n$ an arbitrary non-negative integer). The extension to a non-integer exponent was done in \cite{ZhXi}. 

\item In \cite{ZBL} it was studied the one-dimensional consolidation problem with a threshold gradient. This problem is reduced to a one-phase Stefan problem with a latent heat expressed as $L=\dfrac{\gamma}{\dot{s}(t)}$. That is to say the latent heat depends on the rate of the moving boundary. It must be noted that the case considered in \cite{ZBL} is not properly a Stefan problem because the velocity of the moving boundary disappears, and it has to be treated as a free boundary problem with implicit conditions (\cite{Fa}, \cite{Sc}).

\end{itemize}

Based in the bibliography mentioned above it is quite natural from a mathematical point of view to define a one-phase Stefan problem with a latent heat given by $L=\gamma s^{\beta}(t)\dot{s}^{\delta}(t)$ (with $\gamma$ a given constant and $\beta$ and $\delta$ arbitrary real constants).

It is worth pointing out that this formulation constitutes a  mathematical generalization of the one-phase classical Stefan problem and the problems studied in \cite{VSP}, \cite{ZhXi} and  \cite{ZBL}.

The aim of this paper is to proof in Section 3 the existence and uniqueness of the explicit solution of the problem given by equations (\ref{EcCalor})-(\ref{FreeBound}) in Section 2. Moreover, in Section 4 we will consider some special cases and computational examples for different values of the parameters involved in the problem (\ref{EcCalor})-(\ref{FreeBound}).  The analytical solution that will be obtained in this work will recover in one formula, by  choosing different values for $\beta$ and $\delta$, the solutions obtained in: the classical Stefan problem ($\beta=\delta=0$), the problem considered in  \cite{VSP}( $\beta=1$, $\delta=0$) , the problems solved in \cite{ZBL} and \cite{ZhXi} ( $\beta\in \mathbb{R}^{+}_0$, $\delta=0$) and the problem studied in \cite{ZBL} ( $\beta=0$, $\delta=-1$).

\section{Formulation of the problem.}

$ $

 This paper is intended to study the one-dimensional one-phase Stefan problem for the fusion of a semi-infinite material $x>0$ in which it is  involved a variable latent heat. From a mathematical point of view the problem can be formulated as follows: Find  the free boundary $s=s(t)$ (separation between phases) and the temperature $T=T(x,t)$  in the liquid portion of the material that satisfy  the one-dimensional heat conduction equation:
\be \label{EcCalor}
a^2 T_{xx}(x,t)=T_t(x,t), \qquad 0<x<s(t), 
\ee
 subject to the following boundary condition:
\begin{eqnarray}
T(0,t) &=&t^{\alpha/2}T_0, \label{TempIni} \qquad t>0
\end{eqnarray}
the temperature condition at the interface:
\be 
T(s(t),t)= 0,  \qquad t>0 \label{TempCambioFase}
\ee
the Stefan condition at the interface:
\be
-kT_x(s(t),t)=L(s(t),\dot{s}(t)) \dot{s}(t), \qquad t>0  \label{CondStefan}
\ee
and the initial condition:
\be
s(0)=0. \label{FreeBound}
\ee
Here the parameters $a^2$ (diffusion coefficient) and $k>0$ (thermal conductivity) are known constants.
The phase change temperature is 0 and the imposed temperature at the fixed face $x=0$ is given by $t^{\alpha/2}T_0>0$, where we assume that $\alpha$  is a non-negative real exponent.

The remarkable feature of this problem is related to the condition at the interface given by the Stefan condition (\ref{CondStefan}), where the latent heat by unit of volume will be defined by:

\be
L(s(t),\dot{s}(t))=\gamma s(t)^{\beta} \dot{s}(t)^{\delta}, \label{LatentHeat}
\ee
where $\gamma$ is a given constant, and $\beta$ and $\delta$ are arbitrary real constants.

\section{Explicit solution of the problem.}

In order to solve the problem (\ref{EcCalor})-(\ref{FreeBound}) presented in Section 2  , we use the similarity transformation given by \cite{ZWB} and \cite{ZhXi}: 

\be
T(x,t)=t^{\alpha/2}\varphi\left( \eta\right), \quad \text{with } \quad \eta=\dfrac{x}{2a\sqrt{t}}. \label{Transform-1}
\ee
Computing the derivatives of $T$:
\be
T_{xx}(x,t)=\dfrac{t^{(\alpha/2-1)}}{4a^2} \varphi'' (\eta) \quad \text{and} \quad  T_t(x,t)=\dfrac{\alpha}{2}t^{(\alpha/2-1)}\varphi
(\eta)-t^{(\alpha/2-1)} \varphi'(\eta) \dfrac{\eta}{2},
\ee
we obtain that the temperature given by (\ref{Transform-1}) satisfies the heat equation (\ref{EcCalor}) if and only if $\varphi$ is the solution of the following ordinary differential equation:
\be
\varphi''(\eta)+2\eta \varphi'(\eta)-2\alpha \varphi(\eta) =0, \label{EcuDif-1}
\ee
whose general solution, in this case, can be written as (see the proof in the Appendix A):

\be\label{genSol2}
\varphi(\eta)=c_1 M \left(-\dfrac{\alpha}{2},\dfrac{1}{2},-\eta^2\right)+c_{2}\eta M\left(-\dfrac{\alpha}{2}+\dfrac{1}{2},\dfrac{3}{2},-\eta^2\right),
\ee
where $c_{1}$ and $c_{2}$ are arbitrary  constants. The function $M(a,b,z)$, called Kummer function, is  defined by:

\be 
 M(a,b,z)=\sum\limits_{s=0}^{\infty}\frac{(a)_s}{(b)_s s!}z^s,  \label{M} \\
\ee 
in which b cannot be a non-positive integer, and where $(a)_s$ is the Pochhammer symbol  defined by:
\be 
 (a)_s=a(a+1)(a+2)\dots (a+s-1), \quad \quad (a)_0=1 
\ee
All the properties of Kummer's functions to be used in this paper can be found in \cite{OLBC}.

Therefore  $T(x,t)$ is given by:

\be\label{Temperatura}
T(x,t)=t^{\alpha/2}\left[ c_1 M\left(-\dfrac{\alpha}{2},\dfrac{1}{2},-\eta^2 \right)+c_2 \eta M\left( -\dfrac{\alpha}{2}+\dfrac{1}{2},\dfrac{3}{2}, -\eta^2\right)\right]. 
\ee
where $c_1$ and $c_2$ are constants that must be determined in order that $T(x,t)$ satisfies the conditions (\ref{TempIni})-(\ref{CondStefan}).

From equation (\ref{TempIni}), taking into account that $M\left(-\dfrac{\alpha}{2},\dfrac{1}{2},0 \right)=1$, it is obtained:
\be \label{c1}
c_1=T_0.
\ee

From condition (\ref{TempCambioFase}), there must be $\varphi\left(\dfrac{s(t)}{2a\sqrt{t}} \right)=0$, $\forall t>0$. Then, it can be said that:
\be \label{fronteraLibre}
s(t)=2a\xi \sqrt{t}, 
\ee
where $\xi$ is a positive constant that has to be determined.  Bearing in mind that the free boundary $s(t)$ is defined by (\ref{fronteraLibre}), it can be deduced from 
(\ref{TempCambioFase}) and (\ref{c1}) that:

\be \label{c2}
c_2=\dfrac{-T_0 M\left(-\dfrac{\alpha}{2},\dfrac{1}{2},-\xi^2 \right)}{\xi M\left( -\dfrac{\alpha}{2}+\dfrac{1}{2},\dfrac{3}{2},-\xi^2\right)}.
\ee

Until know, $s(t)$ and $c_2$ are given in function of $\xi$. In order to determine $\xi$,  it will be applied the Stefan condition (\ref{CondStefan}) which has not been considered yet. For that purpose, $T_x(x,t)$ must be calculated:
\be \label{derivTemperatura}
 T_x(x,t)=\dfrac{t^{(\alpha-1)/2}}{a}\left[c_1 \alpha M\left( -\dfrac{\alpha}{2}+1,\dfrac{3}{2},-\eta^2\right)+\dfrac{c_2}{2}M\left( -\dfrac{\alpha}{2}+\dfrac{1}{2},\dfrac{1}{2},-\eta^2\right) \right]
\ee

From the Stefan condition, taking into account (\ref{fronteraLibre}) and (\ref{derivTemperatura}), it is obtained:

\be \label{condStefan-2}
-\dfrac{k t^{(\alpha-1)/2}}{a}\left[c_1 \alpha M\left( -\dfrac{\alpha}{2}+1,\dfrac{3}{2},-\xi^2\right)+\dfrac{c_2}{2}M\left( -\dfrac{\alpha}{2}+\dfrac{1}{2},\dfrac{1}{2},-\xi^2\right) \right]= \gamma 2^{\beta} \xi^{\beta+\delta+1}a^{\beta+\delta+1}t^{(\beta-\delta-1)/2}.
\ee
As $c_1$, $c_2$ and $\xi$ does not depend on $t$, (\ref{condStefan-2}) makes sense if and only if:  $t^{(\alpha-1)/2}=t^{(\beta-\delta-1)/2}$, leading to the relationship:
\be\label{alpha-beta-delta} 
\alpha=\beta-\delta.
\ee

Therefore, assuming that (\ref{alpha-beta-delta}) holds, condition (\ref{CondStefan}) leads to:

\begin{eqnarray}
&&\dfrac{kT_0}{2a\xi M\left( -\dfrac{\alpha}{2}+\dfrac{1}{2},\dfrac{3}{2},-\xi^2\right)} \left[ -2\alpha \xi^2 M\left(-\dfrac{\alpha}{2}+1,\dfrac{3}{2},-\xi^2 \right)M\left( -\dfrac{\alpha}{2}+\dfrac{1}{2},\dfrac{3}{2},-\xi^2\right)+ \right. \nonumber\\ 
&&\qquad \qquad\qquad\left. +M\left( -\dfrac{\alpha}{2},\dfrac{1}{2},-\xi^2\right)M\left(-\dfrac{\alpha}{2}+\dfrac{1}{2},\dfrac{1}{2},-\xi^2 \right) \right]=\gamma 2^{\beta}a^{\beta+\delta+1}\xi^{\beta+\delta+1}. \label{condStefan-3}
\end{eqnarray}
Taking note of the  identity proved in \cite{ZhXi}:
\be 
 e^{-\xi^2}=  -2\alpha \xi^2 M\left( -\dfrac{\alpha}{2}+\dfrac{1}{2},\dfrac{3}{2},-\xi^2\right)M\left(-\dfrac{\alpha}{2}+1,\dfrac{3}{2},-\xi^2 \right) +M\left(-\dfrac{\alpha}{2},\dfrac{1}{2},-\xi^2 \right)M\left( -\dfrac{\alpha}{2}+\dfrac{1}{2},\dfrac{1}{2},-\xi^2\right)
\ee
and the relationship presented in \cite{OLBC}:
\be
M\left( -\dfrac{\alpha}{2}+\dfrac{1}{2},\dfrac{3}{2},-\xi^2\right) =e^{-\xi^2} M\left( \dfrac{\alpha}{2}+1,\dfrac{3}{2},\xi^2\right)
\ee
equation (\ref{condStefan-3}) becomes:
\be 
\dfrac{kT_0}{\gamma a^{\beta+\delta+2}2^{\beta+1}} \dfrac{1}{\xi M\left(\dfrac{\alpha}{2}+1,\dfrac{3}{2},\xi^2 \right)}= \xi^{\beta+\delta+1}.
\ee
That is to say: $\xi$ is a positive solution of the equation:
\be \label{ecuacionXi}
\dfrac{kT_0}{\gamma a^{\beta+\delta+2}2^{\beta+1}} f(z)= z^{\beta+\delta+1}, \qquad z>0,
\ee
where:
\be
f(z)=\dfrac{1}{z M\left(\dfrac{\alpha}{2}+1,\dfrac{3}{2},z^2 \right)}
\ee
Furthermore, knowing from \cite{OLBC} that:
\be 
\dfrac{d}{dz} \left[ zM\left(\dfrac{\alpha}{2}+1,\dfrac{3}{2},z^2 \right)\right]=M\left( \dfrac{\alpha}{2}+1,\dfrac{1}{2},z^2\right)
\ee 
it follows that:
\be
f'(z)=-f^2(z) M\left( \dfrac{\alpha}{2}+1,\dfrac{1}{2},z^2\right).
\ee
In this way, it can be said that the left hand side of equation (\ref{ecuacionXi}) given by $LE(z)=\dfrac{kT_0}{\gamma a^{\beta+\delta+2}2^{\beta+1}} f(z)$ verifies:
\begin{eqnarray}
&&(LE)'(z)=-\dfrac{kT_0}{\gamma a^{\beta+\delta+2}2^{\beta+1}}f^2(x) M\left( \dfrac{\alpha}{2}+1,\dfrac{1}{2},z^2\right)<0, \qquad (\text{since } z>0 \text{ and } \alpha\geq 0) \label{LE-1}\\
&& LE(0)= +\infty, \label{LE-2}\\
&& LE(+\infty)= 0, \label{LE-3}
\end{eqnarray}
meanwhile the right hand side of (\ref{ecuacionXi}) given by $RI(z)=z^{\beta+\delta+1}$ satisfies:
\begin{eqnarray}
&&(RI)'(z)= (\beta+\delta+1)z^{\beta+\delta}>0, \qquad (\text{ if } \beta+\delta+1\geq 0) \label{RI-1}\\
&& RI(0)= 0, \label{RI-2}\\
&& RI(+\infty)= +\infty. \label{RI-3}
\end{eqnarray}
Thus, from (\ref{LE-1})-(\ref{LE-3}) and (\ref{RI-1})-(\ref{RI-3}), one can conclude that equation (\ref{ecuacionXi}) has a unique positive solution $\xi$ provided that $\beta+\delta+1\geq 0$.

It should be mentioned that due to (\ref{alpha-beta-delta}), i.e $\alpha=\beta-\delta$, the fact that $\alpha\geq 0$, and the request that $\beta+\delta+1\geq 0$, the results obtained in this paper are valid only if $\beta\geq \max\left(\delta, -\delta-1 \right)$.

The above arguments can be summarized in the following theorem:

\begin{teo} \label{Teo2.1}
Let $\beta$ and $\delta$ be arbitrary real constants satisfying $\beta\geq \max\left(\delta, -\delta-1 \right)$. Taking $\alpha=\beta-\delta$, there exists a unique solution  of a similarity type for the one-phase Stefan problem (\ref{EcCalor})-(\ref{FreeBound})  given by:
\begingroup
\addtolength{\jot}{0.35em}
\begin{flalign}
&T(x,t)=  t^{\alpha/2}\left[ c_{1} M\left(-\frac{\alpha}{2}, \frac{1}{2},-\eta^2\right)+c_{2} \eta M\left(-\frac{\alpha}{2}+\frac{1}{2},\frac{3}{2},-\eta^2\right)\right], \label{Temp-Teo}\\
&s(t) =2a \xi \sqrt{t}, \label{Frontera-Teo}
\end{flalign}
\endgroup
where  $\eta=\dfrac{x}{2a\sqrt{t}}$ and the constants $c_{1}$ and $c_{2}$ are given by:
\be 
c_{1}=T_0, \qquad  \qquad c_{2}=\dfrac{-T_0 M\left(-\dfrac{\alpha}{2},\dfrac{1}{2},-\xi^2 \right)}{\xi M\left( -\dfrac{\alpha}{2}+\dfrac{1}{2},\dfrac{3}{2},-\xi^2\right)},
\ee
and the dimensionless coefficient $\xi$ is obtained as the unique positive solution of the following equation:

\be  
\dfrac{kT_0}{\gamma a^{\beta+\delta+2}2^{\beta+1}} f(z)= z^{\beta+\delta+1}, \qquad z>0, \label{ecu-teo}
\ee
in which $f$ is the real function defined by:
\be
f(z)=\dfrac{1}{z M\left(\dfrac{\alpha}{2}+1,\dfrac{3}{2},z^2 \right)}, \qquad z>0. \label{f-teo}
\ee
\end{teo}

\section{Special cases and computational examples.}

\par This section is meant to highlight the problems that are generalized in this work by showing that the solutions already reached in the literature can be obtained from the one we present by just choosing the appropriate parameters $\beta$, $\delta$ and thus $\alpha$. For each case it it going to be done a computational example in order to see how the parameter $\xi$, that characterizes the free boundary, varies with respect to $\delta$, for a fixed $\beta$.

\noindent The properties  found in \cite{OLBC} and \cite{ZhXi} will be helpful in the subsequent arguments:
\begin{eqnarray}
M\left(0,b,z \right)&=&1 \label{Prop-1}\\
M\left(a,b,z\right)&=&e^z M(a,b,-z) \label{Prop-2}\\
M\left(-\dfrac{n}{2},\dfrac{1}{2},-z^2 \right)  &=&2^{n-1} \Gamma\left(\dfrac{n}{2}+1 \right)\left[i^n erfc(z)+i^nerfc(-z) \right]\label{Prop-3} \qquad \text{with } n \text{ integer} \\
zM\left(-\dfrac{n}{2}+\dfrac{1}{2},\dfrac{3}{2},-z^2 \right)&=&2^{n-2}\Gamma\left( \dfrac{n}{2}+\dfrac{1}{2}\right)\left[i^n erfc(-z)-i^n erfc(z) \right] \qquad \text{with } n \text{ integer}\label{Prop-4}
\end{eqnarray}
where $i^nerfc(\cdot)$ is the repeated integral of the complementary error function defined by:
\begingroup
\addtolength{\jot}{0.1em}
\begin{align}
& i^0 erfc(z)=erfc(z)=1- erf(z), \qquad erf(z)=\dfrac{2}{\sqrt{\pi}}\int_0^z e^{-u^2}du, \\
& i^n erfc(z)=\int\limits_{z}^{+\infty} i^{n-1}erfc(t)dt 
\end{align}
\endgroup

Let us analyse the explicit solution achieved in each of the following problems:

\vspace{3cm}
\begin{enumerate}

\item \textbf{Classical one-phase Stefan problem}

In the classical Stefan problem the latent heat is given by $L=\gamma$ constant, so the solution in this case can be recovered from the solution given by (\ref{Temp-Teo})-(\ref{f-teo}) by taking $\beta=0, \delta=0$ and thus $\alpha=0$. It must be pointed out that in this case, using the fact that:
\be
M\left(0,\dfrac{1}{2},-\eta^2 \right) =1, \qquad \text{and} \qquad M\left( \dfrac{1}{2},\dfrac{3}{2},-\eta^2\right)=\dfrac{\sqrt{\pi}}{2\eta} erf(\eta),
\ee
we get a temperature with the form $T(x,t)=c_1+c_2 erf(\eta)$, like in the classical literature (\cite{AlSo}, \cite{Ca}, \cite{Cr}, \cite{Gu}, \cite{Lu}, \cite{Ru} and \cite{Ta4}).

\item \textbf{Problem of the shoreline movement in a sedimentary basin \cite{VSP} and \cite{SaTa}}

In this problem it was considered a latent heat that varies linearly with the position, that is to say $L=\gamma s(t)$. Therefore the solution can be obtained from Theorem \ref{Teo2.1} by choosing $\beta=1$, $\delta=0$ and thus $\alpha=1$. Using (\ref{Prop-3}), and the properties: $ierf(z)=\frac{e^{-z^2}}{\sqrt{\pi}}+z erfc(z)$ and $\Gamma\left( \frac{3}{2}\right)=\frac{\sqrt{\pi}}{2}$, the temperature becomes:

\begin{eqnarray}
T(x,t)	  &=&  c_1\left[ \sqrt{t}e^{-\eta^2} +\dfrac{\sqrt{\pi}}{2}x erfc(\eta)  \right] + \dfrac{c_2}{2} x  \label{alpha=1}
\end{eqnarray}
in accordance to the solution shown in \cite{VSP} and \cite{SaTa}.

\item \textbf{Problem with a latent heat defined as a power function of the position solved in \cite{ZWB} and \cite{ZhXi}}

In these papers, the latent heat $L$ is defined as a power function of the position, i.e, $L=\gamma s(t)^{\beta}$ with $\gamma$ constant and $\beta$ a non-negative real exponent. Choosing $\beta\in \mathbb{R}_0^{+}$, $\delta=0$, and then $\alpha=\beta$, the solutions given in \cite{ZWB} and  \cite{ZhXi} are automatically recovered.

\item \textbf{One-dimensional consolidation problem with threshold gradient \cite{ZBL}}

In this work, it is considered a one-dimensional consolidation problem with a threshold gradient which can be transformed into a one-phase Stefan problem with a latent heat that depends on the rate of change of the moving boundary. It is studied the case in which $L=\dfrac{\gamma}{\dot{s}(t)}$.  We remark here that this problem is not a Stefan problem because the velocity of the free boundary does not appear but it is a free boundary problem for the heat equation with implicit free boundary conditions (\cite{Fa}, \cite{Sc}). Fixing $\beta=0$, $\delta=-1$ and so $\alpha=1$, the solution of this problem can be obtained from Theorem \ref{Teo2.1} the temperature can be expressed as (\ref{alpha=1}). In addition, taking into account (\ref{Prop-2}), it is obtained that:

\be 
f(z)= \dfrac{1}{zM\left( \dfrac{3}{2},\dfrac{3}{2},z^2\right)} =\dfrac{1}{z e^{z^2} M\left(0,\dfrac{3}{2},-z^2 \right)}= \dfrac{e^{-z^2}}{z}
\ee 
in agreement to the solution in \cite{ZBL}.

\vspace{0.2cm}

Once it has been compared the solution obtained in this paper with the solutions presented in the literature, we are going to run some computational examples.
In order to solve the Stefan problem (\ref{EcCalor})-(\ref{FreeBound}) it is necessary to solve  equation (\ref{ecu-teo}) which is equivalent to find the unique zero of the following function: 

\be  
H(z)=\dfrac{kT_0}{\gamma a^{\beta+\delta+2}2^{\beta+1}} f(z)- z^{\beta+\delta+1}, \qquad z>0, \label{H-Newton}
\ee
in which $f$ is the function defined by (\ref{f-teo}).

We are going to apply Newton's method  with the iteration formula:

\be 
z_k=z_{k-1}-\dfrac{H(z_{k-1})}{H'(z_{k-1})}.
\ee 
in order to solve (\ref{H-Newton}). For the computational examples it will be considered  the corresponding thermal parameters for the water in liquid state i.e., $k=0.58$ $[W/(m ^{\circ}C)]$ and $a^2=1.39\times 10^{-7}$ $[m^2/s]$. Without loss of generality it will be assumed $\gamma=1$. Newton's Method will be implemented using Matlab software in order to find the unique positive solution of equation (\ref{ecu-teo}). It is worth pointing out that in this programming language, the Kummer function $M(a,b,z)$ is represented by the `hypergeom' command. The stopping criterion to be used here is the boundedness of the   absolute error $\vert z_k-z_{k-1}\vert<10^{-15}$.

Fig.1 shows the variation of $\xi$ (solution of (\ref{H-Newton})) with respect to $\delta$, choosing different values for the coefficient that characterizes the temperature at the fixed face ($T_0=1,5$ or $10$ $[^{\circ}C/s^{\alpha/2}]$) and fixing $\beta=0$.

It is clearly obvious that looking at Fig.1, it is obtained for $\delta=0$ the solution of $\xi$ for the Classical Stefan problem; for $\delta=1$ it is found the solution given by \cite{VSP} and for any other $\delta \in \mathbb{R}_0^{+}$ it is recovered the solution given by \cite{ZWB} and \cite{ZhXi}.

\begin{Figure}
 \centering
 \includegraphics[width=\textwidth]{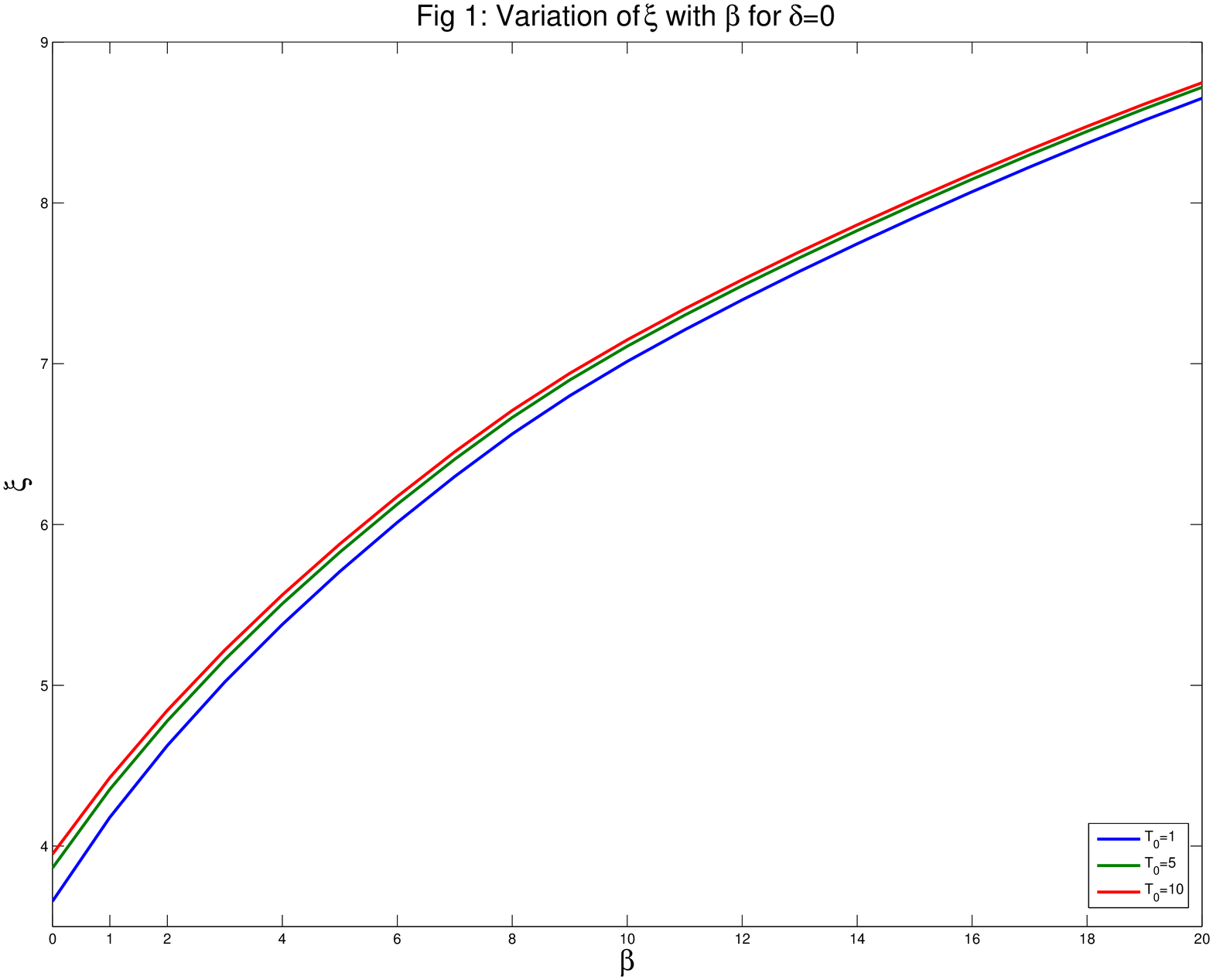}
\end{Figure}

In the same way, Fig. 2, shows the variation of $\xi$ with respect to $\beta$ ($-1\leq \beta\leq 1$), choosing different values for the coefficient that characterizes the temperature at the fixed face ($T_0=1,5$ or $10$ $[^{\circ}C/s^{\alpha/2}]$) and fixing $\delta=0$. The case for $\beta=-1$ corresponds to the solution of the problem analysed in \cite{ZBL}.

\begin{Figure}
\includegraphics[width=\textwidth]{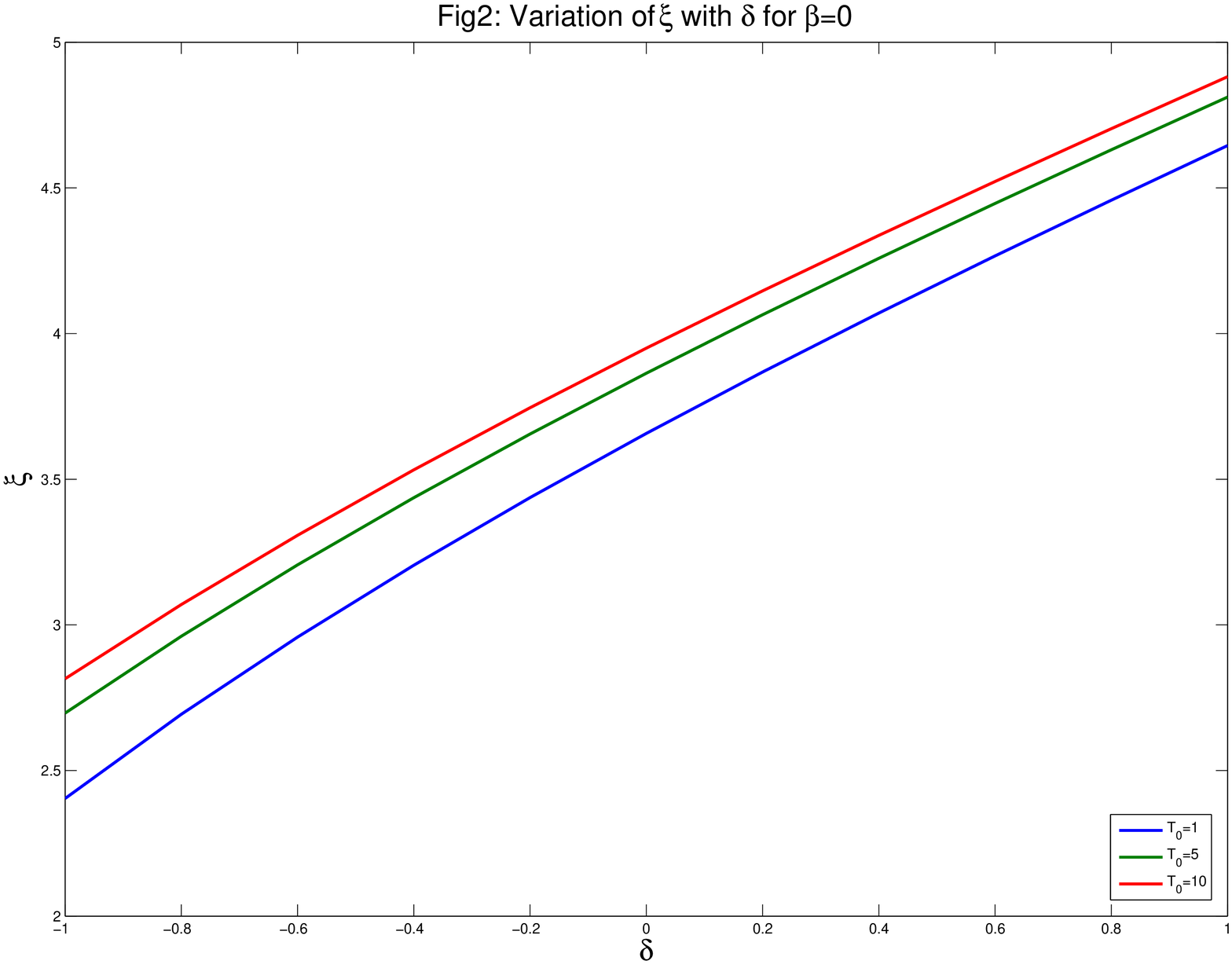}
\end{Figure}

$ $ 

The results obtained  indicate that $\xi$ increases with $\beta$ increasing and $\delta=0$ and the same happens when $\beta=0$ is fixed and $\delta$ varies between -1 and 1.
Moreover, it can be assured that the greater the value of $T_0$, the higher is the value obtained for $\xi$ (parameter that characterizes the free boundary) implying that the fusion process occurs faster.
\end{enumerate}


\section{Conclusions}
$ $

In this work a Stefan problem with a latent heat that depends on the position of the free boundary as well as on its rate of change has been analysed
An explicit solution has been found using the similarity technique and the theory of Kummer functions. This exact solution gathers in one formula the solutions obtained in the previous papers: 
\cite{VSP}, \cite{ZhXi}, \cite{ZWB} and \cite{ZBL}, constituting a generalization of them. Besides, the exact solution is worth finding since it can be used to provide a benchmark for verifying the accuracy of numerical methods that approximate the solution of Stefan problems.

We have also applied Newton's method to the problem (\ref{EcCalor})-(\ref{FreeBound}), in order to estimate the parameter $\xi$ that characterizes the free front numerically. It has been recovered the solutions given in the literature. In addition it was observed that this parameter increases with respect to the parameter $\beta$, fixing $\delta=0$ and vice versa. Also, it can be noted that if the coefficient that characterizes the initial temperature $T_0$ becomes greater, $\xi$ also does, meaning that the phase-change happens quicker, validating mathematically what it seems obvious from the physical point of view.

\appendix
\numberwithin{equation}{section}

\section{Appendix}

It is going to be proved that the general solution of the following ordinary differential equation:
\be
\varphi''(\eta)+2\eta \varphi'(\eta)-2\alpha \varphi(\eta) =0. \label{A-Ecu}
\ee
is given by:
\be\label{A-Sol}
\varphi(\eta)=c_1 M \left(-\dfrac{\alpha}{2},\dfrac{1}{2},-\eta^2\right)+c_{2}\eta M\left(-\dfrac{\alpha}{2}+\dfrac{1}{2},\dfrac{3}{2},-\eta^2\right),
\ee
regardless $\alpha$ is an integer or a non-integer non-negative number, where $c_1$ and $c_2$ are arbitrary constants.

\vspace{0.3cm}
\noindent \textbf{Case $\alpha$  non-negative, non-integer:}
\vspace{0.3cm}

Introducing the new variable $w(\eta)=-\eta^2$ like in \cite{ZhXi} and defining $g(w)=\varphi(\eta(w))$ it is obtained that (\ref{A-Ecu}) is equivalent to the Kummer's differential equation :
\be
wg''(w)+g'(w)\left( \dfrac{1}{2}-w\right) +\dfrac{\alpha}{2}f(w)=0. \label{A-Ecu-2}
\ee
whose general solution according to \cite{OLBC} is given by:
\be
g(w)=\widehat{c_{1}}M \left(-\dfrac{\alpha}{2},\dfrac{1}{2},w\right)+\widehat{c_{2}}U\left(-\dfrac{\alpha}{2},\dfrac{1}{2},w\right),
\ee
where $\widehat{c_1}$ and $\widehat{c_2}$ are arbitrary constants.
Due to the fact that $U$ can be defined as:
\be 
U(a,b,z)=\frac{\Gamma(1-b)}{\Gamma(a-b+1)}M(a,b,z)+\frac{\Gamma(b-1)}{\Gamma(a)} z^{1-b}M(a-b+1,2-b,z) \label{U}.
\ee
it is obtained that the general solution of (\ref{A-Ecu-2}) is given by:
\be 
g(w)=\overline{c_{1}}M \left(-\dfrac{\alpha}{2},\dfrac{1}{2},w\right)+\overline{c_{2}}w^{1/2} M\left(-\dfrac{\alpha}{2}+\dfrac{1}{2},\dfrac{3}{2},w\right),
\ee
where $\overline{c_{1}}$ and $\overline{c_{2}}$ are arbitrary constants, arriving in this way to a $\varphi$ solution of (\ref{A-Ecu}) defined by (\ref{A-Sol}).

\vspace{0.3cm}
\noindent \textbf{Case $\alpha=n$ non-negative integer:}
\vspace{0.3cm}

According to \cite{ZWB}, the general solution of (\ref{A-Ecu}) is given by:
\be
\varphi(\eta)=  \widehat{c_1} i^n erfc(\eta)+\widehat{c_2}i^n erfc(-\eta). \label{A-Sol-3}
\ee
where $i^n erfc(\cdot)$ is the family of the repeated integrals of the complementary error function.

Let $c_1$ and $c_2$ arbitrary constants. Taking: 
\begin{eqnarray}
\widehat{c_1}&=&c_1 2^{n-1}\Gamma\left(\dfrac{n}{2}+1 \right)-c_2 2^{n-2} \Gamma\left(\dfrac{n}{2}+\dfrac{1}{2} \right) \\
\widehat{c_2}&=&c_1 2^{n-1}\Gamma\left(\dfrac{n}{2}+1 \right)+c_2 2^{n-2} \Gamma\left(\dfrac{n}{2}+\dfrac{1}{2} \right)
\end{eqnarray}
in (\ref{A-Sol-3}) it leads to:
\begin{eqnarray}
\varphi(\eta)&=&\left[c_1 2^{n-1}\Gamma\left(\dfrac{n}{2}+1 \right)-c_2 2^{n-2} \Gamma\left(\dfrac{n}{2}+\dfrac{1}{2} \right) \right] i^n erfc(\eta)+ \nonumber\\
&&+  \left[c_1 2^{n-1}\Gamma\left(\dfrac{n}{2}+1 \right)+ c_2 2^{n-2} \Gamma\left(\dfrac{n}{2}+\dfrac{1}{2} \right) \right] i^n erfc(-\eta), \nonumber\\
&=& c_1 2^{n-1} \Gamma\left(\dfrac{n}{2}+1 \right)\left[ i^n erfc(\eta)+i^n erfc(-\eta) \right]+ \nonumber \\
&&+c_2 2^{n-2} \Gamma\left(\dfrac{n}{2}+\dfrac{1}{2} \right)\left[ i^n erfc(-\eta)-i^n erfc(\eta) \right], \nonumber\\
&=& c_1 M\left( -\dfrac{n}{2},\dfrac{1}{2},-\eta^2\right)+c_2\eta M\left(-\dfrac{n}{2}+\dfrac{1}{2}, \dfrac{3}{2},-\eta^2\right).
\end{eqnarray} 
arriving to a solution of (\ref{A-Ecu}) given by a $\varphi$ defined as (\ref{A-Sol}), using the properties stated in (\ref{Prop-3})-(\ref{Prop-4}).

\section*{Aknowledgements}

The present work has been partially sponsored by the Project PIP No 0275 from CONICET-UA, Rosario, Argentina, and Grant AFOSR-SOARD FA9550-14-1-0122.

\small{
}

\end{document}